\documentclass[12pt]{article}
\usepackage{amsfonts}
\usepackage{amsmath}
\usepackage{amssymb}
\usepackage{theorem}

\newtheorem{theorem}{Theorem}
\newtheorem{proposition}[theorem]{Proposition}
\newtheorem{lemma}[theorem]{Lemma} 
\newtheorem{corollary}[theorem]{Corollary}
\theorembodyfont{\rmfamily}
\newtheorem{example}[theorem]{Example}

\newtheorem{remark}[theorem]{Remark}
\newcommand{\Aut}{{\mathrm{Aut}}}

\newcommand{\Inj}{{\mathrm{Inj}}}
\newcommand{\maxx}{{\mathrm{max}}}
\newcommand{\Out}{{\mathrm{Out}}}
\newcommand{\semi}{{\rtimes}}
\newcommand{\act}{{\cdot}}
\newcommand{\calf}{{\mathcal F}}

\title{Realising fusion systems}
\author{Ian J. Leary\thanks{Partially supported by NSF grant
DMS-0505471}
\and Radu Stancu}

\date{\today}

\newenvironment{proof}[1][]{\begin{trivlist} \item[\hskip\labelsep
\emph{Proof#1.}]}{\foorp \end{trivlist}}
\newcommand{\foorp}{{\unskip\nobreak\hfil\penalty50
 \hskip1em\vadjust{}\nobreak\hfil $\Box$
 \parfillskip=0pt \finalhyphendemerits=0 \par}}

\begin{document} 

\maketitle

\begin{abstract} 
We show that every fusion system on a $p$-group $S$ is equal to the 
fusion system associated to a discrete group $G$ with the property 
that every $p$-subgroup of $G$ is conjugate to a subgroup of $S$.  
\end{abstract}

\section{Introduction} 
Let $p$ be a prime number.  By a $p$-group we shall mean a 
finite group whose order is a power of $p$.  A 
\emph{fusion system} on a $p$-group $S$ is a category $\calf$ 
whose objects are the subgroups of $S$, and whose morphisms 
are injective group homomorphisms, subject to certain axioms. 
The notion of a fusion system is intended to axiomatize 
the $p$-local structure of a discrete group $G\geq S$ 
in which every $p$-subgroup is conjugate to a subgroup of $S$.  
Every such $G$ gives rise to a fusion system $\calf_S(G)$ 
on $S$, and we say that $G$ realises $\calf$ if $\calf_S(G)=\calf$.  

The notion of a \emph{saturated fusion system} is intended to 
axiomatize the $p$-local structure of a finite group in which 
$S$ is a Sylow $p$-subgroup.  It is known that there are 
saturated fusion systems $\calf$ which are not realised by 
any finite group $G$, although showing that this is the case 
is very delicate.  In the case when $p=2$, the only known 
examples are certain systems discovered by Ron Solomon \cite{djb,lo,sol}.  

In contrast, we show that every fusion system on any 
$p$-group $S$ is realised by some discrete group $G\geq S$ 
in which every maximal $p$-subgroup is conjugate to $S$.  
The groups $G$ that are used in our proofs are constructed 
as graphs of finite groups.  In particular each of our 
groups $G$ contains a free subgroup of finite index.  In 
an appendix we give a brief account of those parts of the 
theory of graphs of groups that we use.  

While preparing this paper, we learned that Geoff Robinson 
has proved a similar, but not identical result \cite{rob}.  
Since~\cite{rob} was already submitted when we started to 
write this paper, we have taken it upon ourselves to compare
and contrast the two results.  
Robinson's construction realises a large class of fusion systems, 
including all saturated fusion systems, but does not realise 
all fusion systems.  The groups that Robinson constructs are 
iterated free products with amalgamation, whereas the groups 
that we construct are iterated HNN extensions.  In both 
cases the groups may be viewed as graphs of finite groups.  

We state and outline the proof of a version of Robinson's 
theorem, along the lines of the proof of our main result.  
We also give examples of fusion systems that cannot be 
realised by Robinson's method, we give examples of 
non-saturated fusion systems that are realised by Robinson's 
method, and we prove an analogue of Cayley's theorem for 
fusion systems.  

The work in this paper grew from the authors' participation 
in the Banff conference `Homotopy theory and group actions' 
and from a VIGRE reading seminar at Ohio State which studied 
the Aschbacher-Chermak approach to the Solomon fusion
systems~\cite{aschcherm}.  The authors thank Andy Chermak 
and Geoff Robinson for showing them early versions of 
\cite{aschcherm}~and~\cite{rob}.  

\section{Definitions and results} 
Let $p$ be a prime, and let $G$ be a discrete group.  The
$p$-Frobenius category $\Phi_p(G)$ of the group $G$ is a 
category whose objects are the $p$-subgroups of $G$.  
If $P$ and $Q$ are $p$-subgroups of $G$, or equivalently 
objects of $\Phi_p(G)$, the morphisms from $P$ to $Q$ are 
the group homomorphisms $f:P\rightarrow Q$ that are equal 
to conjugation by some element of $G$.  Thus $f:P\rightarrow Q$ 
is in $\Phi_p(G)$ if and only if there exists $g\in G$ with 
$f(u)= g^{-1}ug$ for all $u\in P$.  (Note that the element 
$g$ is \emph{not} part of the morphism.  If $g'=zg$ for 
some element $z$ in the centralizer of $P$, then $g$ and $g'$ 
define the same morphism.)  

Now suppose that $S$ is a $p$-subgroup of $G$ that is 
maximal, and further suppose that every $p$-subgroup of $G$ 
is conjugate to a subgroup of $S$.  In this case, every object 
of $\Phi_p(G)$ is isomorphic within the category 
$\Phi_p(G)$ to a subgroup of $S$.  It follows that the 
full subcategory $\calf_S(G)$ with objects the subgroups of $S$ 
is equivalent to $\Phi_p(G)$.  This example motivates Puig's definition
of a \emph{fusion system on $S$} \cite{puig}.  A fusion system 
on a $p$-group $S$ is a category $\calf$.  The objects of $\calf$ 
are the subgroups of $S$, and the morphisms from $P$ to $Q$ form 
a subset of the set $\Inj(P,Q)$ of injective group homomorphisms
from $P$ to $Q$.  These are subject to the following axioms:  
\begin{enumerate} 
\item For any $s\in S$, and any $P, Q\leq S$ with $s^{-1}Ps\leq Q$, 
the morphism $\phi:P\rightarrow Q$ defined by 
$\phi:u\mapsto s^{-1}us$ is in $\calf$; 
\item If $f:P\rightarrow Q$ is in $\calf$, with $R=f(P)\leq Q$, 
then so are $f:P\rightarrow R$ and $f^{-1}:R\rightarrow P$.  
\end{enumerate}  
It is easily checked that these axioms are satisfied in the case 
when $\calf=\calf_S(G)$ as defined above.  Note that the 
first axiom could be rewritten as the statement $\calf_S(S)\subseteq 
\calf$.  

\begin{remark} 
Fusion systems arise in other ways.  For example, if $H$ is any 
group and $S$ is any $p$-subgroup of $H$, then the full 
subcategory of $\Phi_p(H)$ with objects the subgroups of $S$ 
is a fusion system on $S$.  Another source of fusion systems on a
$p$-group $S$ is the Brauer category of a $p$-block~$b$~\cite{alpbrou,linck}.  
Here $H$ is a finite group, $S$ is the defect group of the 
$p$-block $b$, and the morphisms in the category are those 
conjugations by elements of $H$ that preserve some extra structure 
associated to $b$.  In the case when $b$ is the principal block, 
$S$ is the Sylow $p$-subgroup of $H$ and this fusion system is 
just $\calf_S(H)$.  One corollary of our Theorem~\ref{main} 
is that every such fusion system is realised by some group $G$.  
\end{remark}

There is a fusion system $\calf_{S}^{\,\maxx}$ on $S$, in which the 
morphisms from $P$ to $Q$ consists of all injective group
homomorphisms from $P$ to $Q$.  Any fusion system on $S$ is a 
subcategory of $\calf_{S}^{\,\maxx}$, and the intersection of a family 
of fusion systems on $S$ is itself a fusion system.  
If $\Phi=\{\phi_1,\ldots,\phi_r\}$ is a collection of morphisms in
$\calf_{S}^{\,\maxx}$, where $\phi_i:P_i\rightarrow Q_i$, the 
fusion system generated by $\Phi$ is defined to be the 
smallest fusion system that contains each~$\phi_i$.  

\begin{theorem}
\label{main} 
Suppose that $\calf$ is the fusion system on $S$ generated by
$\Phi=\{\phi_1,\ldots,\phi_r\}$.  Let $T$ be a free group with 
free generators $t_1,\ldots,t_r$, and define $G$ as the quotient 
of the free product $S*T$ by the relations $t_i^{-1}ut_i=\phi_i(u)$ 
for all $i$ and for all $u\in P_i$.  Then $S$ embeds as a subgroup 
of $G$, every $p$-subgroup of $G$ is conjugate to a subgroup of $S$, 
and $\calf_S(G)=\calf$.  Moreover, 
every finite subgroup of $G$ is conjugate to a subgroup of $S$, and 
$G$ has a free normal subgroup of index dividing $|S|!$.  
\end{theorem} 

If $f:S'\rightarrow S$ is an injective group homomorphism between 
$p$-groups, and $\calf'$ is a fusion system on $S'$, then 
there is a functor $f_*$ from $\calf'$ to $\calf_{S}^{\,\maxx}$, which 
sends $P'\leq S'$ to $f(P')$ and $\phi':P'\rightarrow Q'$ to 
$$f\circ \phi'\circ f^{-1}: f(P')\rightarrow f(Q').$$

\begin{theorem}[Robinson~\cite{rob}] 
\label{two} 
Suppose that $\calf$ is the fusion system on $S$ generated by 
the images $(f_i)_*(\calf_{S_i}(G_i))$ for injective group homomorphisms 
$f_i:S'_i\rightarrow S$ for $1\leq i\leq r$, where $G_i$ is a 
finite group with $S'_i$ as a Sylow $p$-subgroup.  Define $G$ as 
the quotient of the free product $S*G_1*\cdots *G_r$ by the 
relations $s=f_i(s)$ for all $i$ and for all $s\in S'_i$.  Then 
$S$ embeds as a subgroup of $G$, 
every $p$-subgroup of $G$ is conjugate to a subgroup of $S$, and 
$\calf_S(G)=\calf$.  Moreover, every finite subgroup of $G$ is 
conjugate to a subgroup of one of the $G_i$, or to a subgroup of $S$, 
and $G$ has a free normal subgroup of index dividing $N!$, where $N$ 
is the least common multiple of $|S|$ and the $|G_i|$.  
\end{theorem} 

\begin{remark} 
The above theorem can be obtained from theorem~1 of~\cite{rob} by
induction.  The main result of~\cite{rob} is theorem~2, which is 
similar to the above statement except that extra conditions are 
put on the $G_i$.  
\end{remark} 

\begin{theorem} \label{cayley} 
Let $\Sigma$ denote the group of all permutations of the elements of a 
$p$-group $S$, and identify $S$ with a subgroup of $\Sigma$ via the 
Cayley embedding.  Every fusion system on $S$ is equal to a 
subcategory of the Frobenius category $\Phi_p(\Sigma)$ of $\Sigma$.  
\end{theorem} 

\section{Saturated fusion systems} 

In this section we present the definition of a saturated fusion 
system, due to Puig~\cite{puig}, although we shall describe an 
equivalent definition due to Broto, Levi and Oliver~\cite{blo}.  
There are two additional axioms as well as the axioms for a 
fusion system.  These axioms necessitate some preliminary definitions.  

As usual, if $G$ is a group and $H$ is a subgroup of $G$, we write 
$C_G(H)$ for the centralizer of $H$ in $G$ and $N_G(H)$ for the 
normalizer of $H$ in $G$.  

Suppose that $\calf$ is a 
fusion system on $S$.  Say that $P\leq S$ is fully $\calf$-centralized 
if 
$$|C_S(P)|\geq |C_S(P')|$$ 
for every $P'$ which is isomorphic to $P$ as an object of $\calf$.  
Suppose that 
$\calf= \calf_S(G)$ for some discrete group $G$ in which every 
$p$-subgroup is conjugate to a subgroup of $S$.  In this case, if 
$P$ is fully $\calf$-centralized, one sees that $C_S(P)$ is a 
$p$-subgroup of $C_G(P)$ of maximal order.  

Similarly, say that $P$ is fully $\calf$-normalized if 
$$|N_S(P)|\geq |N_S(P')|$$ 
for every $P'$ which is isomorphic to $P$ as an object of $\calf$.  
If $\calf=\calf_S(G)$ as above and $P$ is fully $\calf$-normalized, 
one sees that $N_S(P)$ is a $p$-subgroup of $N_G(P)$ of maximal 
order.  

Now suppose that $\calf=\calf_S(G)$ for some finite group $G$, 
and that $P\leq S$ is fully $\calf$-normalized.  In this case, 
$N_S(P)$ must be a Sylow $p$-subgroup of the finite group 
$N_G(P)$.  Moreover, $C_G(P)\cap N_S(P)= C_S(P)$ must be a 
Sylow $p$-subgroup of $C_G(P)$, and $\Aut_S(P)=N_S(P)/C_S(P)$ 
must be a Sylow $p$-subgroup of $\Aut_\calf(P)=N_G(P)/C_G(P)$.  
This gives the first of two extra axioms for a saturated fusion system: 

\begin{enumerate}
\setcounter{enumi}{2}
\item 
If $P$ is fully $\calf$-normalized, then $P$ is also fully
$\calf$-centralized, and $\Aut_S(P)$ is a Sylow $p$-subgroup 
of $\Aut_\calf(P)$.  
\end{enumerate}

Next, suppose that $\calf=\calf_S(G)$ for some finite group $G$ 
and that $f:P\rightarrow Q\leq S$ is an isomorphism in $\calf$ 
such that $Q$ is fully $\calf$-centralized.  This implies that 
$C_S(Q)$ is a Sylow $p$-subgroup of $C_G(Q)$.  Pick an element 
$h\in G$ so that $f$ is equal to conjugation by $h$, i.e., so 
that $f(u)=c_h(u)=h^{-1}uh$ for all $u\in P$.  The image 
$c_h(C_S(P))$ is a $p$-subgroup of $C_G(c_h(P))=C_G(Q)$, and 
so there exists $h'\in C_G(Q)$ so that $c_{h'}\circ c_h(C_S(P))
\leq C_S(Q)$.  Since $c_{h'}$ acts as the identity on $Q$, 
if we define $k=hh'$, we see that $c_k$ extends $f$ and 
$c_k(C_S(P))\leq C_S(Q)$.  

The map $c_k$ clearly extends to a map from 
$N_f=N_S(P)\cap c_k^{-1}(N_S(Q))$ to $N_S(Q)$.  But since 
$C_S(P)$ is a subgroup of $c_k^{-1}(N_S(Q))$, we may rewrite this as 
$$N_f = \{ g\in N_S(P) \colon c_k\circ c_g\circ c_k^{-1}\in
\Aut_S(Q)\} $$
$$\phantom{N_f} = \{ g\in N_S(P) \colon f\circ c_g\circ f^{-1}\in
\Aut_S(Q)\},$$ 
which does not depend on choice of $k$.  This leads to the 
second extra axiom: 

\begin{enumerate}
\setcounter{enumi}{3}
\item 
If $f:P\rightarrow Q$ is an isomorphism in $\calf$ and $Q$ is 
fully $\calf$-centralized, then $f$ extends in $\calf$ to a map 
from $N_f$ to $N_S(Q)$, where 
$$N_f = \{ g\in N_S(P) \colon f\circ c_g\circ f^{-1}\in
\Aut_S(Q)\}. $$ 
\end{enumerate}

\begin{remark} It has been shown~\cite{kessarstancu} 
that the axioms for a saturated fusion system can be simplified to: 
\begin{enumerate} 
\item[3$'$.] $\Aut_S(S)$ is a Sylow $p$-subgroup of $\Aut_\calf(S)$. 
\item[4$'$.] If $f:P\rightarrow Q$ is an isomorphism in $\calf$ and $Q$ is 
fully $\calf$-normalized, then $f$ extends in $\calf$ to a map 
from $N_f$ to $N_S(Q)$, where $N_f$ is as defined in axiom~4.
\end{enumerate} 
\end{remark}

\begin{remark}\label{abelian} 
In the case when $S$ is abelian, axioms 3~and~4 simplify.  
In this case, every subgroup of $S$ is fully $\calf$-centralized 
and fully $\calf$-normalized for any fusion system $\calf$, and for 
any $f\in \calf$, $N_f= S$.  Hence  
a fusion system $\calf$ on an abelian $p$-subgroup $S$ is saturated if 
and only if $\Aut_\calf(S)$ is a $p'$-group and every morphism
$f:P\rightarrow S$ in $\calf$ extends to an automorphism of $S$.  
\end{remark} 

\begin{remark} 
As mentioned in the introduction, there are saturated fusion systems 
which are not realised by any finite group.  One source of 
saturated fusion systems is the fusion systems associated to 
$p$-blocks of finite groups~\cite{alpbrou,linck}.  The question 
of whether every such fusion system can be realised by a finite 
group is a long-standing open problem.  
\end{remark} 

\section{Examples} 

Let $E$ be an elementary abelian $p$-group of rank at least three, i.e.,
a direct product of at least three copies of the cyclic group of order
$p$.  Let $A=\Aut(E)$ be the full group of automorphisms of $E$, which 
is of course isomorphic to a general linear group over the field of
$p$ elements.  Let $B$ be a subgroup of $A$ of order a power of $p$, 
and let $C$ be a non-trivial subgroup of $A$ of order coprime to~$p$.  
Note that $A$ is generated by its subgroups of order coprime to~$p$.  

Each of $A$, $B$ and $C$ may be viewed as a collection of morphisms in 
the fusion system $\calf_{E}^{\,\maxx}$.  For $X=A$, $B$ or $C$, let
$\calf_E(X)$ denote the fusion system generated by all the morphisms 
in $X$.  

\begin{example} \label{egone} 
The fusion system $\calf_E(C)$ 
is saturated, and is equal to the fusion system 
$\calf_E(G)$, where $G$ is the semi-direct product $G=E\semi C$.  
\end{example} 

\begin{example} \label{egtwo} 
The fusion system $\calf_E(A)$ is not 
saturated, since in $\calf_E(A)$ the 
automorphism group of the object $E$ does not have $E/Z(E)$ as a 
Sylow $p$-subgroup.  However, $\calf_E(A)$ can be realised by 
the procedure of Theorem~\ref{two}.  Let $C_1,\ldots,C_r$ be 
$p'$-subgroups of $A$ that together generate $A$.  If we put 
$G_i=E\semi C_i$ with $f_i$ the identity map of $E$, then the  
fusion system generated by all of the $(f_i)_*(\calf_E(G_i))$ 
is equal to $\calf_E(A)$.  
\end{example} 

\begin{example} \label{egthree} 
The fusion system $\calf_E(B)$ cannot be realised by the procedure 
used in Theorem~\ref{two}.  For suppose that $G_1,\ldots,G_r$ are 
finite groups with Sylow $p$-subgroups $E_1,\ldots,E_r$, each of 
which is isomorphic to a subgroup of $E$, and suppose that 
$\calf_E(B)$ is generated by the fusion systems $(f_i)_*\calf_{E_i}(G_i)$.  
Those $G_i$ for which $f_i:E_i\rightarrow E$ is not an isomorphism 
do not contribute any 
morphisms to $\Aut_\calf(E)$.  If $f_i:E_i\rightarrow E$ is an 
isomorphism, then either $\Aut_{G_i}(E_i)$ contains non-identity 
elements of $p'$ order, implying that $\calf\neq \calf_E(B)$, or 
$E_i$ is central in $G_i$ and $G_i$ does not contribute any 
morphisms to $\Aut_\calf(E)$.  
\end{example} 

Next we consider some examples of fusion systems $\calf$ on an 
abelian $p$-group $E$ in which $\Aut_\calf(E)$ is a $p'$-group, 
but for which some isomorphisms between proper subgroups of $E$ 
do not extend to elements of $\Aut_\calf(E)$.  

\begin{example}\label{egfour} 
Let $F$ and $F'$ be distinct order $p$ subgroups of $E$, and let 
$\phi:F\rightarrow F'$ be an isomorphism.  Let $\calf_E(\phi)$ be 
the fusion system generated by $\phi$.  Every morphism in
$\calf_E(\phi)$ is equal to either an inclusion map or the 
composite of either $\phi$ or $\phi^{-1}$ with an inclusion map.  
In particular, in $\calf_E(\phi)$, the automorphism group of each 
object $E'\leq E$ is trivial.  The fusion system $\calf_E(\phi)$ 
cannot be realised by the procedure of Theorem~\ref{two}, as will 
be explained below.  

In view of Remark~\ref{abelian}, $\calf_E(\phi)$ is not a 
saturated fusion system, since the 
morphism $\phi:F\rightarrow F'$ does not extend to an automorphism 
in $\calf_E(\phi)$ of the group $E$.  

Now suppose that $\calf$ is a fusion system on $E$ generated by 
the images $(f_i)_*\calf_{E_i}(G_i)$ of some fusion systems for 
finite groups.  If $\phi:F\rightarrow F'$ is a morphism in 
$\calf$, then there exists 
$i$ so that $F, F'\leq f_i(E_i)$ and $\phi\in
(f_i)_*\calf_{E_i}(G_i)$.  But then (by the same argument as 
used above) there is a morphism $\tilde\phi:f_i(E_i)\rightarrow 
f_i(E_i)$ extending $\phi:F\rightarrow F'$.  
Thus $\calf$ cannot be equal to the fusion system $\calf_E(\phi)$, 
since this fusion system contains no such $\tilde\phi$.  
\end{example} 

\begin{example} \label{egfive} 
Let $F$ be a proper subgroup of $E$, and suppose that $D$ is a 
non-trivial $p'$-group of automorphisms of $F$.  Let $F\semi D$ 
denote the semi-direct product of $F$ and $D$, let $G$ be 
the free product with amalgamation $G=E*_F(F\semi D)$, and let 
$\calf$ be the fusion system $\calf_E(G)$.  From this definition 
one sees that $\calf$ can be obtained by the procedure of
Theorem~\ref{two}. On the other hand, since $\Aut_\calf(E)$ 
is trivial, one sees that the non-trivial automorphisms 
of $F$ do not extend to automorphisms of $E$, and hence $\calf$ 
is not saturated.  
\end{example} 

As remarked earlier, Robinson does not consider all fusion systems
that can be built by the procedure of Theorem~\ref{two}, but only
those fusion systems that he calls Alperin fusion systems~\cite{rob}.
With the notation of Theorem~\ref{two}, a fusion system is Alperin if
the following conditions hold:
\begin{enumerate} 
\item Inside each $G_i$ there is a subgroup $E_i$ which is the 
largest normal $p$-subgroup of $G_i$, and the centralizer of this 
subgroup is as small as possible, in the sense that 
$C_{G_i}(E_i)= Z(E_i)$; 
\item The quotient $G_i/E_i$ is isomorphic to $\Out_\calf(E_i)
:=\Aut_\calf(E_i)/\Aut_{E_i}(E_i)$; 
\item Inside $S$, the image of the subgroup $S'_i$ (the Sylow
$p$-subgroup of $G_i$ which is to be identified with a subgroup 
of $S$) is equal to the 
normalizer of the image of $E_i$, i.e., $f_i(S'_i)= N_S(f_i(E_i))$. 
\end{enumerate} 

In terms of this definition, the content of Alperin's fusion theorem 
with some later embellishments~\cite{alperin,goldschmidt} 
is that the fusion system for any finite group is Alperin.  
Robinson remarks~\cite{rob} 
that work of Broto, Castellana, Grodal, Levi and Oliver 
implies that every saturated fusion system is Alperin~\cite{bcglo}.  
It is easy to see that a fusion system on an abelian $p$-group is 
Alperin if and only if it is saturated.  We finish this section 
by giving an 
example of a fusion system that is Alperin but not saturated.  

\begin{example} \label{egsix} 
Let $p$ be an odd prime, let $A= (C_p)^3$, and let $B$ be a subgroup 
of $\Aut(A)$ of order $p$ such that $A$ is indecomposable as a
$B$-module.  (Equivalently, the action of a generator for $B$ on $A$ 
should be a single Jordan block.)  Let $S$ be the semi-direct product 
$S=A\semi B$.  The centre $Z$ of $S$ has order $p$.  Let $E = Z\times
B \leq S$, a subgroup isomorphic to $C_p\times C_p$.  It is readily 
seen that $C_S(E)=E$ and that $P=N_S(E)$ is isomorphic to a
semi-direct product $(C_p)^2\semi C_p$, the unique non-abelian group 
of order $p^3$ and exponent $p$.  Let $G_1$ be the semi-direct product 
$G_1=E\semi \Aut(E)$.  Since the Sylow $p$-subgroups of $\Aut(E)$ are 
cyclic of order $p$, there is an isomorphism between $P$ and a 
Sylow $p$-subgroup of $G_1$ that extends the inclusion of $E$.  

By construction, the fusion system $\calf$ for the free product 
with amalgamation $S*_PG_1$ is Alperin 
in the sense of Robinson~\cite{rob}, but this fusion system is not 
saturated.  For example, there are non-identity self-maps of $Z$
inside $\calf$, and if $\calf$ were saturated, any self-map of $Z$ 
inside $\calf$ would extend to a self-map of $S$.  But in $\calf$, 
$S$ has only inner automorphisms, and these restrict to $Z$ as the 
identity.   
\end{example} 

\section{Proofs} 

\begin{proof} (of Theorem~\ref{cayley}.)  As in the statement, let 
$\Sigma$ be the group of all permutations of $S$, and identify $S$
with a subgroup of $\Sigma$.  Let $P$ and $Q$ be subgroups of $S\leq
\Sigma$, and let $\phi:P\rightarrow Q$ be any injective group
homomorphism.  It suffices to show that there is some $\sigma\in
\Sigma$ such that for all $u\in P$, $\sigma^{-1} u\sigma=\phi(u)$.
Let $\Omega$ denote the group $S$ viewed as a set with a left
$S$-action.  There are two ways to view $\Omega$ as a set with a left
$P$-action, via $P\leq S$ and via $\phi:P\rightarrow Q\leq S$.  Denote
these two $P$-sets by $\Omega$ and $^\phi\Omega$ respectively.  Each
of $\Omega$ and $^\phi\Omega$ is isomorphic as a $P$-set to the
disjoint union of $|S:P|$ copies of $P$.  In particular, there is an
isomorphism of $P$-sets $\sigma:{}^\phi\Omega\rightarrow\Omega$.
Viewing $\sigma$ as an element of $\Sigma$, one has that $\sigma
\phi(u)\omega=u\sigma\omega$ for all $u\in P$ and $\omega\in \Omega$.
Hence $\sigma^{-1} u\sigma =\phi(u)$ for all $u$ as required.
\end{proof} 

\begin{remark} A version of Theorem~\ref{cayley} appeared
  in~\cite{lsy}, although fusion systems were not mentioned there.  
\end{remark} 

Before proving Theorem~\ref{main} we give a result concerning
extending group homomorphisms, and two corollaries, one of 
which will be used in the proof.  

\begin{lemma} \label{homext} 
Let $S$ and $G$ be as in the statement of
Theorem~\ref{main}, let $j:S\rightarrow G$ be the natural map from 
$S$ to $G$, let $H$ be a group and let $f:S\rightarrow H$ be a group
homomorphism.  There is a group homomorphism $\tilde f:G\rightarrow 
H$ with $f=\tilde f\circ j$ if and only if for each $i$, the homomorphisms
$f:P_i\rightarrow H$ and $f\circ \phi_i:P_i\rightarrow H$ differ by an inner
automorphism of $H$.   
\end{lemma} 

\begin{proof} 
Given a homomorphism $\tilde f$ as in the statement, one has that 
for each $i$ and for each $u\in P_i$, 
$f\phi_i(u) = h_i^{-1} f(u) h_i$, where $h_i = \tilde f(t_i)$.  
For the converse, suppose that there exists, for each $i$, an element 
$h_i$ satisfying the equation 
$f\phi_i(u)= h_i^{-1} f(u) h_i$ for all $u\in P_i$.  In this case 
one may define $\tilde f$ on the generators of $G$ by $\tilde
f(s)=f(s)$ for all $s\in S$ and $\tilde f(t_i)=h_i$.  
\end{proof} 

\begin{corollary} \label{maptosymm} 
With notation as in the statement of Theorem~\ref{main}, there is a 
homomorphism from $G$ to $\Sigma$, the group of all permutations of 
the set $S$, extending the Cayley representation of $S$.  
\end{corollary} 

\begin{proof} 
The argument used in the proof of Theorem~\ref{cayley} shows that 
the conditions of Lemma~\ref{homext} hold.
\end{proof} 

\begin{remark} Corollary~\ref{maptosymm} gives an alternative way to 
prove Corollary~\ref{vertinj}, at least in the special case of a
rose-shaped graph.  
\end{remark} 

\begin{corollary} 
With notation as in the statement of Theorem~\ref{main}, a complex 
representation of $S$ with character $\chi$ 
extends to a complex representation of $G$ if and only if for each 
$i$ and for each $u\in P_i$, $\chi(u)= \chi(\phi_i(u))$.  
\end{corollary} 

\begin{remark} Of course, a representation of $S$ will extend to $G$ 
in many different ways if it extends at all.  
\end{remark}

\begin{proof} (of Theorem~\ref{main}.) 
As in Appendix~\ref{pressec}, one sees that the group $G$ presented 
in the statement is the fundamental group of a graph of groups with one
vertex group, $S$, and one edge group $P_i$ for each $\phi_i$, $1\le i\le r$.   
From Corollary~\ref{vertinj} it follows that $S$ is a subgroup of 
$G$.  From Corollary~\ref{finsubgps}, it follows that any finite
subgroup of $G$, and in particular any $p$-subgroup of $G$, is 
conjugate to a subgroup of $S$.  By Theorem~\ref{tree}, there is 
a cellular action of $G$ on a tree $T$, with one orbit of vertices and 
$r$ orbits of edges.  By suitable choice of orbit representatives,
we may choose a vertex $v$ whose stabilizer is $S$, and edges 
$e_1,\ldots,e_r$ so that the stabilizer of $e_i$ is $P_i$, and 
so that the initial vertex of $e_i$ is $v$ while the final vertex 
is $t_i\act v$.  

Since every $p$-subgroup of $G$ is conjugate to a subgroup of $S$, 
there is a fusion system $\calf_S(G)$ associated to $G$.  By
construction $\calf_S(G)$ contains each $\phi_i$, which corresponds 
to conjugation by $t_i$.  

Conversely, suppose that $g\in G$ has the 
property that $g^{-1}Pg\leq Q$ for some subgroups $P,Q$ of $S$.  It 
suffices to show that conjugation by $g$, as a map from $P$ to $Q$, 
is equal to a composite of (restrictions of) the maps $\phi_j$ and 
their inverses with conjugation maps by elements of $S$.  

Consider the action of $P$ on the tree $T$.  By hypothesis, the 
action of $P$ fixes both the vertex $v$ and the vertex $g\act v$.  
Since $T$ is a tree, $P$ must fix all the 
vertices and edges on the unique shortest path from $v$ to $g\act v$.  
Let this path have 
length $n$.  Define $g_0=1_G$, $g_n=g$, and for $1\leq i\leq n-1$,
choose $g_i\in G$ so that $g_0\act v,g_1\act v,\ldots,g_n\act v$ is the shortest 
path in $T$ from $v$ to $g\act v$.   For each $i$, $P$ is contained in 
the stabilizer of the vertex $g_i\act v$, and so $P\leq g_iSg_i^{-1}$, 
or equivalently $g_i^{-1}Pg_i \leq S$.  

The edge joining $g_i\act v$ and 
$g_{i+1}\act v$ is an edge of the form $g_i\act e_j$ or $g_{i+1}\act e_j$ for 
some $j$ depending on $i$.  Consider the two cases separately, first 
supposing that the edge is of the form $g_i\act e_j$.  In this case 
it follows that 
$P\leq g_iP_ig_i^{-1}$, since $P$ stabilizes the edge $g_i\act e_j$.  Also 
one sees that $g_{i+1}\act v=g_it_j\act v$, and hence $g_{i+1}^{-1}g_it_j\in
S$.  Hence conjugation by $g_i^{-1}g_{i+1}$, viewed as a map from 
$g_i^{-1}Pg_i$ to $g_{i+1}^{-1}Pg_{i+1}$ is equal to the composite 
of the map $\phi_j$ (restricted to $g_i^{-1}Pg_i \leq P_i$) followed
by conjugation by an element of $S$.  

The other case is similar.  Here it follows that $P\leq
g_{i+1}P_{i+1}g_{i+1}^{-1}$, and one has that $g_i\act v= g_{i+1}t_j\act v$, 
from which $g_i^{-1}g_{i+1}t_j=s\in S$.  In this case conjugation by 
$g_i^{-1}g_{i+1}$, as a map from $g_i^{-1}Pg_i$ to
$g_{i+1}^{-1}Pg_{i+1}$, is equal to the composite map given by 
conjugation by $s$ followed by the map $\phi_j^{-1}$ (restricted 
to $s^{-1}g_i^{-1}Pg_is\leq \phi_j(P_{i+1})$).  

Thus conjugation by $g=g_n$ as a map from $P$ to $Q$ can be expressed 
as a composite of maps inside the fusion system generated by the 
$\phi_i$, and so $\calf_S(G)$ is equal to this fusion system.  

It remains to show that the group $G$ contains a free normal subgroup 
of index at most $|S|!$.  Let $\Sigma$ denote the symmetric group on 
the set $S$.  By Corollary~\ref{maptosymm}, there is a homomorphism 
$G\rightarrow \Sigma$ which extends the natural injection
$S\rightarrow \Sigma$.  By Corollary~\ref{virtfree}, the kernel of 
this homomorphism is a free normal subgroup of $G$, and its index is 
a factor of $|\Sigma|=|S|!$.  
\end{proof} 

\begin{proof} (of Theorem~\ref{two}---sketch.)  In this case, the group 
$G$ is the fundamental group of a star-shaped graph of groups, with 
one central vertex labelled $S$ and $r$ outer vertices labelled
$G_1,\ldots, G_r$.  The edge from $G_i$ to $S$ is labelled by the 
group $S'_i$. By Theorem~\ref{tree}, there is 
a cellular action of $G$ on a tree $T$, with $r+1$ orbit of vertices and 
$r$ orbits of edges. We may choose orbit representatives
$v_0,v_1,\ldots,v_r$ of vertices and $e_1,\ldots,e_r$ of edges so 
that the stabilizer of $v_0$ is $S$, and for $1\leq i\leq r$, the 
stabilizer of $v_i$ is $G_i$ (resp.\ of $e_i$ is $S'_i$).  Moreover, 
we may assume that $e_i$ has initial vertex $v_i$ and terminal vertex 
$v_0$.  

In this case, one sees that any finite subgroup of $G$ is conjugate 
to either a subgroup of $S$ or to a subgroup of $G_i$ for some $i$.  
Since $S'_i$ is a Sylow $p$-subgroup of $G_i$, it follows that any 
$p$-subgroup of $G$ is conjugate to a subgroup of $S$ as required.  

As in the previous proof, it is clear that the fusion system
$\calf_S(G)$ contains the image of each $\calf_{S'_i}(G_i)$, but an
argument is needed to show that these images generate $\calf_S(G)$.
Given $g\in G$ and $P,Q\leq S$ so that $g^{-1}Pg\leq Q$, one argues 
that the action of $P$ fixes the vertices $v_0$ and $g\act v_0$ in the 
tree $T$, and hence fixes the shortest path (necessarily of even
length, say $2n$) that joins these vertices.  

Let $g_0=1_G$, $g_{2n}=g$, and pick group elements so that the vertices 
on the shortest path from $v_0$ to $g\act v_0$ are: 
$$g_0\act v_0\,,\,\, g_1\act v_{j(1)}\,,\,\,g_2\act v_0\,,\,\,
g_3\act v_{j(2)}\,,\,\ldots\,,\,g_{2n-1}\act v_{j(n)}\,,\,\,g_{2n}\act v_0$$ 
for some function $j:\{1,\dots,n\}\rightarrow \{1,\ldots,r\}$.  
If $i$ is even, then $g_i^{-1}Pg_i\leq S$, and if $i$ is odd 
then $g_i^{-1}Pg_i\leq G_{j(({i+1})/{2})}\,$.  Since $P$ stabilizes each 
edge, one sees that $P\leq g_i^{-1}S_kg_i$, where $S_k$ denotes 
the image of $S'_k$ inside $S$, and $k=j(({i+1})/{2})$ if $i$ is odd
and $k=j({i}/{2})$ if $i$ is even.  In particular, each $g_i^{-1}Pg_i$ is
a subgroup of $S$.  

One may show that in the case when $i$ 
is odd, $g_i^{-1}g_{i+1}\in G_{j(({i+1})/{2})}$ and that in the case 
when $i$ is even, $g_i^{-1}g_{i+1}\in S$.  Thus the map from 
$g_i^{-1}Pg_i$ to $g_{i+1}^{-1}Pg_{i+1}$ given by conjugation by 
$g_i^{-1}g_{i+1}$ is a map inside the fusion system generated by 
the images of the $\calf_{S'_i}(G_i)$, and conjugation by $g=g_{2n}$
as a map from $P$ to $Q\leq S$ is expressed as a composite of maps 
of the required form.  

Finally, if $\Omega$ is a finite set so that $|\Omega|$ is divisible
by $|S|$ and by each $|G_i|$, one may define free actions of $S$ and 
each $G_i$ on $\Omega$ which give rise to the same (free) action of 
$S_i=f_i(S'_i)$.  This gives rise to a group homomorphism from 
$G$ to $\Sigma$, the symmetric group on $\Omega$, whose kernel is 
free by Corollary~\ref{virtfree}.  

\end{proof}

\section{Appendix: graphs of groups} 

In this section we give proofs of those results about graphs of groups 
that we use.  Our treatment of graphs of groups follows that of 
Scott and Wall~\cite{scottwall}.  

For the purposes of this paper, a graph $\Gamma$ consists of two sets, 
the vertices $V$ and the directed edges $E$, together with two 
functions $\iota,\tau: E\rightarrow V$.  For $e\in E$, $\iota(e)$ 
is called the initial vertex of $e$ and $\tau(e)$ is the terminal
vertex of $e$.  Multiple edges and loops are allowed in this 
definition.  $\Gamma$ is connected if the only equivalence relation 
on $V$ that contains all pairs $(\iota(e),\tau(e))$ is the relation 
with just one class.  

A graph $\Gamma$ may be viewed as a category, with 
objects the disjoint union of $V$ and $E$ and two non-identity 
morphisms with domain $e$ for each $e\in E$, one morphism 
$e\rightarrow \iota(e)$ and one morphism $e\rightarrow \tau(e)$.  

A graph of groups is a connected graph $\Gamma$ together with 
groups $G_v$, $G_e$ for each vertex and edge, and injective group 
homomorphisms $f_{e,\iota}:G_e\rightarrow G_{\iota(e)}$ 
and $f_{e,\tau}:G_e\rightarrow 
G_{\tau(e)}$ for each edge $e$.  If $\Gamma$ is viewed as a 
category, this is just a functor from $\Gamma$ to 
the category of groups and injective group homomorphisms.  
Without loss of generality, one may assume that each map 
$f_{e,\iota}:G_e\rightarrow G_{\iota(e)}$ is the inclusion of a
subgroup.  

\subsection{The fundamental group of a graph of groups} 

For a topologist, and arguably for anybody, the easiest way to 
define the fundamental group of a graph of groups is via the 
notion of a graph of spaces.  

A graph of spaces is a connected graph $\Gamma$ together with 
topological spaces $X_v$, $X_e$ for each vertex and edge, and 
continuous maps $f_{e,\iota}: X_e\rightarrow X_{\iota(e)}$ and 
$f_{e,\tau}:X_e\rightarrow X_{\tau(e)}$.  
A graph of spaces is just a functor from the category $\Gamma$ to 
the category of topological spaces and continuous functions.  
A graph of based spaces is defined similarly: each $X_e$ and $X_v$ 
is equipped with a base point, and the maps must preserve base
points.  Let $I$ denote the closed unit interval $[0,1]$.  
The total space of a graph of spaces is the space $X$ 
made from the disjoint union 
$$\coprod_{v\in V} X_v 
\,\,\,{\textstyle\coprod}\,\,\, 
\coprod_{e\in E}X_e\times I$$ 
by identifying $(x,0)\in X_e\times I$ with $f_{e,\iota}(x)\in
X_{\iota(e)}$ and identifying $(x,1)\in X_e\times I$ with
$f_{e,\tau}(x) \in X_{\tau(e)}$.  As an example, consider the 
graph of spaces in which each $X_e$ and $X_v$ is a single point.  
For this graph of spaces the total space is the usual topological
realization of the graph as a 1-dimensional CW-complex.  
The reader who is familiar with the homotopy colimit construction 
will note that if one views a graph of spaces as a functor $X_{(-)}$
on the category $\Gamma$, then the total space $X$ is naturally 
homeomorphic to the homotopy 
colimit of the functor $X_{(-)}$, or in symbols, 
$X={\rm hocolim}_\Gamma X_{(-)}$.

Given a graph of groups, one may define a graph of connected based 
spaces by taking classifying spaces as the spaces $X_e$ and $X_v$: 
$$X_e=BG_e=K(G_e,1)\qquad X_v=BG_v=K(G_v,1).$$ 
For the continuous map $f_{e,\iota}:X_e\rightarrow X_{\iota(e)}$ 
(resp.\ $f_{e,\tau}:X_e\rightarrow X_{\tau(e)}$) one may 
take any continuous map that induces the given map $G_e\rightarrow 
G_{\iota(e)}$ (resp.\ $G_e\rightarrow G_{\tau(e)}$) on fundamental 
groups.  Define a total space $X$ as the realization of this graph 
of spaces.  

For discrete groups $K$ and $H$, the space $BK$ is unique 
up to based homotopy, and homotopy class of based maps from $BK$ to 
$BH$ are in bijective correspondence with group homomorphisms from 
$K$ to $H$.  It follows that the homotopy type of the space $X$ 
defined above depends only on the graph of groups, rather than on 
the particular choices of classifying spaces and maps between them.  
The fundamental group $G$ of the graph of groups can now be defined as 
the fundamental group of $X$.  This describes the fundamental group 
of the graph of groups up to isomorphism.  The inclusion of each 
$X_v$ in $X$ defines a conjugacy class of homomorphism 
$G_v\rightarrow G$ (which will be shown to be injective, below).  
For many purposes one wants a more precise description of $G$, 
together with a single choice of homomorphism $G_v\rightarrow G$.  
This can be done by choosing a basepoint for the space $X$, and 
for each $v$, a path in $X$ from the basepoint for $X$ to the 
basepoint for $X_v\subseteq X$.  

\subsection{Presentations for graphs of groups} \label{pressec} 

We shall only consider presentations for graphs of groups where 
the underlying graph is either a `rose' or a `star'.  By a rose 
we mean a graph with only one vertex, so that every edge has 
the same initial and terminal vertices.  By a star we mean a 
connected graph with $n+1$ vertices and $n$ edges, for some $n>0$, 
with one central vertex, such that all the edges have this vertex 
as their terminal vertex and so that every other vertex is the 
initial vertex of exactly one edge.  

Suppose one is given a $p$-group $S$, subgroups $P_i, Q_i\leq S$, 
and injective group
homomorphisms $\phi_i:P_i\rightarrow Q_i$ for $1\leq i \leq r$, as in
the statement of Theorem~\ref{main}.  Use this data to make a
rose-shaped graph of groups with $r$ edges.  Let $S$ be the 
vertex group, let $P_i$ be the $i$th edge group, with the 
inclusion map $P_i\leq S$ (resp.\ the composite $\phi_i:P_i\rightarrow
Q_i\leq S$) as the $i$th initial (resp.\ terminal) homomorphism.  
There is a model for $BP_i$ having just one 0-cell and one 1-cell 
for each element of $P_i$.  Take a model for $BS$ having just one 
0-cell and take this 0-cell as the base point.  To make a CW-complex 
of the homotopy type of the total space of the graph of groups, it 
suffices to add to $BS$ one 1-cell $t_i$ for each $i$ (with both 
ends at the unique 0-cell), one 2-cell $D_{i,u}$ for $1\leq i\leq r$ 
and for each $u\in P_i$, and higher dimensional cells (which will not 
affect the fundamental group).  The attaching map for the 2-cell 
$D_{i,u}$ spells out the word $u\,t_i\,\phi_i(u)\,t_i^{-1}$, and so 
the presentation coming from this CW-structure is the presentation 
given in the statement of Theorem~\ref{main}.  

Next suppose that one is given a $p$-group $S$, groups $G_i$ 
for $1\leq i\leq r$ with Sylow $p$-subgroups $S_i$, and injective 
group homomorphisms $f_i:S_i\rightarrow S$, i.e., the data found 
in the statement of Theorem~\ref{two}.  In this case, define a 
star of groups with central vertex group $P$, other vertex groups 
$G_1,\ldots,G_r$, and edge groups $S_1,\ldots,S_r$.  The map of 
each edge group into its initial vertex group is the inclusion 
$S_i\rightarrow G_i$, and the map of each edge group into its 
terminal vertex group is $f_i:S_i\rightarrow S$.  An argument 
similar to that given in the previous paragraph shows that the 
fundamental group of this graph of groups has the presentation 
given in the statement of Theorem~\ref{two}.  Note that here one 
can make a space homotopy equivalent to the total space of the 
graph of spaces by starting from the one-point union of $BS$ 
and the $BG_i$, without adding any extra 1-cells.  This is 
reflected in the fact that the vertex groups generate the 
fundamental group of the graph of groups.  

\subsection{Properties of graphs of groups} 

\begin{proposition} Let $G$ be the fundamental group of a graph of
  groups based on a graph $\Gamma$.  Every subgroup $H\leq G$ is 
itself the fundamental group of a graph of groups, indexed by a 
graph $\Delta$ equipped with a map $f:\Delta\rightarrow \Gamma$ 
which does not collapse any edges.  For each $v$~and $e\in\Delta$, 
the group $H_v$ (resp.\ $H_e$) is a subgroup of $G_{f(v)}$ (resp.\ 
$G_{f(e)}$).  
\end{proposition} 

\begin{proof}  Use the bijection between connected covering spaces 
of a connected CW-complex (with a choice of base point) and subgroups
of its fundamental group.  Let $X$ be the total space of the graph 
of spaces used in the definition of $G$, so that there is a covering
space of $X$ whose fundamental group is $H$.  Any connected covering 
space of $X$ can be expressed as the total space of a graph of spaces 
indexed by some $\Delta$ as in the statement.  This gives an 
expression for the fundamental group of any connected covering 
space of $X$ as the fundamental group of a graph of groups as 
claimed.  
\end{proof} 

\begin{theorem}\label{class} 
Let $X$ be the total space of the graph of spaces 
used in the definition of the fundamental group $G$ of a graph of 
groups.  The universal covering space of $X$ is contractible, 
and hence $X$ is homotopy equivalent to $BG$.  
\end{theorem} 

\begin{proof}  
We shall build a space $Y$, in such a way that it is clear that 
$Y$ is contractible, and that $Y$ is a covering space of $X$.  

For $v$ a vertex, define the subspace $X'_v$ of $X$ by 
$$X'_v= X_v \cup \bigcup_{\iota(e)=v} X_e\times [0,0.5) 
\cup \bigcup_{\tau(e)=v} X_e\times (0.5,1].$$ 
Similarly, define for $e$ an edge, $X'_e= X_e\times (0,1)$.  
The inclusions $X_v\rightarrow X'_v$ and $X_e\cong X_e\times \{0.5\}
\rightarrow X'_e$ are homotopy equivalences, and it may be useful
to think of $X'_v$ as a nice open neighbourhood of $X_v$ in $X$.  
Let $Y_v$, $Y'_v$, $Y_e$, and $Y'_e$ be the universal covering spaces
of $X_v$, $X'_v$, $X_e$ and $X'_e$ respectively.  Each $Y'_v$ (resp.\  
$Y'_e$) is contractible since it is the universal covering space of the 
classifying space $BG_v$ (resp.\ $BG_e$).  

The definition of the 
space $X'_v$ lifts to a description of the space $Y'_v$.  The
complement $Y'_v-Y_v$ is identified with a collection of disjoint 
copies of $Y_e\times (0,0.5)$, and $Y_e\times (0.5,1)$, for different 
edges $e$.  There are copies of $Y_e\times (0,0.5)$ if and only if 
$\iota(e)=v$.  In this case the copies are in bijective correspondence 
with the cosets of $f_{e,\iota}(G_e)$ in $G_v$.  Similarly, there are 
copies of $Y_e\times (0.5,1)$ for each $e$ with $\tau(e)=v$, and these 
copies are indexed by cosets of $f_{e,\tau}(G_e)$ in $G_v$.  

By induction, we shall construct a sequence 
$Y_0\subseteq Y_1\subseteq Y_2 \cdots$ of
spaces so that: each $Y_n$ is contractible; there is a map
$\pi:Y_n\rightarrow X$ which is locally a covering map except at 
some points of $X$; for any $x\in X$ and any $n\geq 0$, at least 
one of $\pi:Y_n\rightarrow X$ and $\pi:Y_{n+1}\rightarrow X$ 
is locally a covering map at $x$.  
 
Pick a vertex $v$ of the graph $\Gamma$, and define $Y_0$ to be the 
space $Y'_v$.  Define a map $\pi:Y'_v\rightarrow X$ as the composite 
of the map $Y'_v\rightarrow X'_v$ and the inclusion $X'_v\subseteq
X$.  As remarked earlier, $Y'_v-Y_v$ consists of lots of subspaces of the 
form $Y_e\times (0,0.5)$ for $\iota(e)=v$ and lots of subspaces of 
the form $Y_e\times (0.5,1)$ for $\tau(e)=v$.  Define $Y_1$ by 
attaching to each such subspace a copy of $Y'_e$.  The map
$\pi:Y_0\rightarrow X$ extends uniquely to $\pi:Y_1\rightarrow X$ 
by insisting that on each newly-added $Y'_e$ subspace, $\pi$ is 
equal to the composite map $Y'_e\rightarrow X'_e\subseteq X$.  
From the construction of $Y_1$, it is apparent that $Y_1$ is 
contractible.  

In constructing $Y_1$, one attached to $Y_0$ many spaces of the 
form $Y'_e$, by identifying one end of $Y'_e$ with part of $Y_0$. 
For each copy of $Y'_e$ that was attached via its initial end, 
take a copy of $Y'_{\tau(e)}$, and attach this at the other end 
of $Y'_e$.  Similarly, for each copy of $Y'_e$ that was attached 
to $Y_0$ by its terminal end, take a copy of $Y'_{\iota(e)}$ and 
attach this at the other end of $Y'_e$.  This defines a space 
$Y_2$, which is clearly contractible, and the map $\pi$ extends 
uniquely to a map $Y_2\rightarrow X$ which agrees with the 
covering map $Y'_e\rightarrow X'_e$ or $Y'_v\rightarrow X'_v$ 
on each such subspace.  

Now suppose that $n$ is even, and that $Y_n$ has been constructed 
from $Y_{n-1}$ by attaching subspaces $Y'_v$ in such a way that 
the intersection of $Y_{n-1}$ and each new $Y'_v$ is equal to 
one of the components of $Y'_v-Y_v$.  Furthermore, suppose that 
the map $\pi$ on each new $Y'_v$ is equal to the map $Y'_v\rightarrow 
X'_v\subseteq X$.  Form $Y_{n+1}$ by attaching 
a copy of $Y'_e$ to each other component of $Y'_v-Y_v$ for each 
of the copies of $Y'_v$.  Extend the map $\pi$ as before.  

In the case when $n$ is odd, suppose that $Y_n$ has been constructed 
from $Y_{n-1}$ by attaching subspaces $Y'_e$ in such a way that the 
intersection of $Y_{n-1}$ and each new $Y'_e$ is equal to one of 
the two components of $Y'_e-Y_e\times\{0.5\}$.  Suppose also that 
the map $\pi$ on each of the new $Y'_e$ is equal to the map 
$Y'_e\rightarrow X'_e\subseteq X$.  Form $Y_{n+1}$ by attaching a
copy of $Y'_v$ to the other component of each $Y'_e-Y_e\times\{0.5\}$, 
where $v$ is either $\iota(e)$ or $\tau(e)$ depending which component 
of $Y'_e-Y_e\times\{0.5\}$ was used.  Extend the map $\pi$ in the 
same way as before.  

By construction, each $Y_n$ is contractible, and comes equipped with 
a map $\pi:Y_n\rightarrow X$.  If $n$ is even, this map is locally a 
covering except possibly at points of $X$ contained in the union of
the images of the $X_v$.  If $n$ is odd, this map is a covering except 
possibly at point of $X$ contained in the union of the images of 
the $X_e\times\{0.5\}$.  Now define $Y$ by $Y=\bigcup_n Y_n$.  This 
space $Y$ is contractible, and the map $\pi:Y\rightarrow X$ is a
covering map, since it is locally a covering map at every point of
$X$.  It follows that $Y$ is the universal covering space of $X$.  
Since the universal covering space of $X$ has been shown to be 
contractible, it follows that $X$ is a model for $BG$.  
\end{proof}

\begin{remark} 
The above proof relies on the fact that the edge groups map 
injectively to the vertex groups.  
\end{remark}

\begin{corollary}\label{vertinj} 
Each vertex group $G_v$ maps injectively into the fundamental group of
a graph of groups.  
\end{corollary} 

\begin{proof} 
Given a vertex $v$, construct the universal covering space as in 
the proof of Theorem~\ref{class}, with $Y_0=Y'_v$.  The group of 
all deck transformations of $Y$ is naturally isomorphic to $G$, 
the fundamental group of $X$.  Under this isomorphism, the subgroup 
of those deck transformations that preserve $Y_0$ is identified 
with $G_v$.  
\end{proof} 

\begin{remark} There is also an algebraic proof that each $G_v$ 
embeds in $G$.  In the case when the graph is a rose, this argument 
is given in Corollary~\ref{maptosymm}.  
\end{remark} 

\subsection{The action on a tree} 

Say that an action of a group on a tree is cellular if no element 
of the group exchanges the ends of any edge.  

\begin{theorem} \label{tree} 
Let $G$ be the fundamental group of a graph of 
groups indexed by the graph $\Gamma$.  There is a tree $T$ with  
a cellular $G$-action and an isomorphism $f:T/G\cong \Gamma$.  If $\tilde x$ 
is either a vertex or edge of $T$, and $x=f(\tilde x)$ is the 
image of $G\act \tilde x$ under $f$, then the stabilizer of $\tilde x$ 
is conjugate to $G_x$.  
\end{theorem} 

\begin{proof} 
Let $X$ be the total space of the graph of spaces used in defining 
$G$.  As remarked earlier, the underlying topological space of 
the graph $\Gamma$ can be identified with the total space of 
the constant graph of 1-point spaces indexed by $\Gamma$.  
The unique map from each $X_v$ and $X_e$ to a point induces a 
map from $X$ to $\Gamma$.  

Now let $Y$ be the universal covering 
space of $X$, as constructed in the proof of Theorem~\ref{class}.  
This $Y$ can be viewed as a graph of spaces over some graph $\Delta$,
with vertex spaces copies of the spaces $Y_v$ and edges spaces 
copies of the spaces $Y_e$.  The group $G$ acts on $Y$ in such a 
way that the setwise stabilizer of each copy of $Y_v$ is a 
conjugate of $G_v$, and 
similarly the setwise stabilizer of each copy of $Y_e\times (0,1)$ 
is a conjugate of $G_e$.  
Define $T$ to be the total space of the graph of 1-point spaces 
over the graph $\Delta$.  By construction, $T$ is a graph 
equipped with a $G$-action, an equivariant map $\phi:Y\rightarrow T$, 
and an isomorphism $f:T/G\rightarrow \Gamma$.  To check that $T$ is a 
tree, let $T_n=\phi(Y_n)$.  As in the proof of Theorem~\ref{class}, 
one shows inductively that $T_n$ is contractible, and $T=\bigcup_n
T_n$.  
\end{proof}

\begin{lemma} 
Any cellular action of a finite group $H$ on a 
tree $T$ fixes a vertex. 
\end{lemma} 

\begin{proof} 
Take any vertex $t\in T$, 
and define a finite subtree $T'$ to be the union of all the 
shortest paths between elements of the orbit $H\act t$.  If $T'$ 
is not itself fixed by $H$, remove an $H$-orbit of `leaves' 
(i.e., vertices of valency one) from $T'$, and continue this 
process until a subtree fixed by $H$ is all that remains.  
\end{proof} 

\begin{corollary}\label{finsubgps} 
Every finite subgroup of the fundamental group of a 
graph of groups is conjugate to a subgroup of a vertex group. 
\end{corollary} 

\begin{proof} 
Let $G$ be the fundamental group of the graph of groups and let $T$ be
the corresponding tree.  If $H$ is a finite subgroup of $G$ then $H$ 
fixes some vertex of $T$.  The stabilizer of each vertex of $T$ is a 
conjugate of one of the vertex groups $G_v$.  
\end{proof} 

\begin{corollary} \label{virtfree} 
Let $H$ be a subgroup of a graph of groups whose intersection with 
each conjugate of each vertex group is trivial.  Then $H$ is a free 
group.  
\end{corollary}

\begin{proof}  
The hypotheses imply that $H$ acts freely on the tree $T$, and so 
the quotient space $T/H$ is a 1-dimensional classifying space for
$H$.  
\end{proof}

\leftline{\bf Authors' addresses:}

\smallskip
Ian J. Leary: {\tt leary@math.ohio-state.edu} 

\smallskip
Radu Stancu: {\tt stancu@math.ohio-state.edu} 

\smallskip
Department of Mathematics, 

The Ohio State University, 

231 West 18th Avenue, 

Columbus, 

Ohio 43210

\end{document}